\documentclass{article}
\usepackage{amsthm}
\usepackage{amsfonts}
\usepackage{cite}
\usepackage{amsmath, amscd}
\newtheorem{theorem}{Theorem}[section]

\newtheorem{proposition}[theorem]{Proposition}
\newtheorem{lemma}[theorem]{Lemma}

\begin{document}
\author{Jeremy Berquist}
\title{Rational and Semi-Rational Singularities}
\maketitle

\noindent
\linebreak
\textbf{Abstract.  }  It is a conjecture of Koll\'ar that a variety $X$ with rational singularities in some open subvariety $U$ has a rationalification; that is, a proper, birational morphism $f: Y \rightarrow X$ such that $Y$ has rational singularities, and which is an isomorphism over $U$.  Whether this is true is already unknown in the case of a (normal) threefold with rational singularities along a curve $C$ except at a single point $p \in C$.  There is an analogous conjecture for demi-normal varieties $X$, where we must insist that $Y$ has only semi-rational singularities.  Our main result is that if a stronger version of Koll\'ar's conjecture is true for rational (normal) singularities, then the analogous conjecture is true for Gorenstein semi-rational (non-normal) singularities.  We illustrate this first for surfaces.   As the same procedure does not carry over directly to higher dimensions, we must use a slightly different proof for that case.  We include conditions under which the Cohen-Macaulay property (slightly weaker than that of having rational singularities) is preserved when passing between a demi-normal variety and its normalization.  It is this condition that is automatic for surfaces and which makes a proof slightly different in higher dimensions.  In fact, it is not the case that the normalization of a Cohen-Macaulay variety is Cohen-Macaulay (or vice-versa, as is already apparent with surfaces).  These auxiliary results are needed for the proof of the general result, but they are interesting in their own right.  The proof of the general statement is at the end of the paper. 

\noindent
\linebreak
\textbf{Acknowledgement.  }  The author would like to thank Karl Schwede and S\'andor Kov\'acs for many helpful suggestions.  

\newpage
\tableofcontents
\addcontentsline{}{}{}
\newpage
\begin{section}{Introduction}  
A normal variety $X$ is said to have rational singularities if there exists a resolution of singularities $f: Y \rightarrow X$ such that $R^if_*\mathcal{O}_Y = 0$ for $i > 0$.  If this is true for one such resolution, then it can be proved that the same holds for any resolution.  Equivalently, $X$ is Cohen-Macaulay and $f_*\omega_Y = \omega_X.$  See \cite{KM98}, 5.10.  Likewise, a demi-normal variety (that is, a variety with properties $G_1, S_2$, and semi-normality) is said to have semi-rational singularities if the analogous statements hold for a semi-resolution $f : Y \rightarrow X$.  We call $f$ a semi-resolution when $Y$ is semi-smooth, $f$ is proper and birational, and no component of the conductor of $Y$ is $f$-exceptional.

Koll\'ar asked whether a rationalification exists for normal varieties with specified rational singularities.  Stated differently for projective varieties, the problem is whether, for a projective closure $V \hookrightarrow Y_0$ of a quasi-projective variety $V$ with rational singularities, there should exist a proper, birational morphism $Y \rightarrow Y_0$ such that $Y$ has only rational singularities and which is an isomorphism over $V$.  There is an analogous problem for semi-rational quasi-projective varieties, and it might be asked whether a semi-rationalification exists, knowing that rationalifications exist.

In this paper, we show that this is true, up to two assumptions.  First, we assume a somewhat stronger version of Koll\'ar's conjecture, motivated by the definition of rational pairs.  In particular, suppose that $V$ is quasi-projective and has rational (normal) singularities.  Suppose also that $C$ is a given reduced divisor on $V$, and which is Cohen-Macaulay. Let $X$ be a compactification of $V$.  Then there is a projective, birational morphism $Y \rightarrow X$, which is an isomorphism over $V$, $Y$ has only rational singularities, and the closure of $C$ in $Y$ is Cohen-Macaulay.  This last statement is assumed to be true in the present paper.

Now, suppose that $U$ is demi-normal and quasi-projective.  Suppose that $U$ has Gorenstein semi-rational singularities.  Then, assuming the above form of Koll\'ar's conjecture, $U$ has a projective closure with semi-rational singularities everywhere.  At present, we do not know if the Gorenstein hypothesis is needed; likewise, whether the stronger version of Koll\'ar's conjecture for pairs is necessary.  Future research may uncover a result that uses neither of these hypotheses. 

The problem in two dimensions turns out to have a slightly different proof than the problem in arbitrary dimensions.  One reason is that the normalization of a demi-normal surface, which is by definition $S_2$ hence Cohen-Macaulay, is automatically Cohen-Macaulay.  In higher dimensions, there exists a Macaulayfication, as was shown by Kawasaki \cite{K00}.  However, it is not true that the normalization of a Cohen-Macaulay variety is again Cohen-Macaulay.  In fact, one can generically project a normal, non-Cohen-Macaulay projective variety to a hypersurface, noting that hypersurfaces are always Cohen-Macaulay (even Gorenstein).  Thus there must be some way to obtain a normal, Cohen-Macaulay projective closure of a normal, Cohen-Macaulay quasi-projective variety.  This would only partly solve the problem in higher dimensions, because it must also be checked that $f_*\omega_Y = \omega_X$ for a resolution of singularities $f: Y \rightarrow X$, provided we can ensure that $X$ is Cohen-Macaulay.  

Our approach to the problem in arbitrary dimensions is slightly different than the surface case.  In fact, one cannot first find a Macaulayfication and then check the condition on dualizing sheaves.  Kawasaki's Macaulayfication does not preserve the Cohen-Macaulay locus.  Thus, as in resolution of singularities, surfaces seem to welcome their own proofs that do not work for other varieties.  The idea of first finding a Cohen-Macaulay projective closure and then arranging for the condition on sheaves to hold does not seem to work with higher-dimensional varieties.  

We are working only with projective versions of rationalification and semi-rationalification.  However, the results here hold also for quasi-projective varieties. Our approach using projective closures is a simplification of the problem, which could equally well be stated for quasi-projective varieties.

The paper is organized as follows.  In the first section, we show that the existence of rationalifications for normal surfaces implies the existence of semi-rationalification for demi-normal surfaces.  The main idea is that a Gorenstein, demi-normal quasi-projective surface is semi-rational if and only if its normalization is rational.  In the second section, we present two results that indicate when a demi-normal variety and its normalization are both Cohen-Macaulay (hence that the main result is true).  In the final section, we give our proof of the main theorem:  if $U$ is Gorenstein and has semi-rational singularities, then it has a projective closure $X$ with only semi-rational singularities.  Our approach is similar to the surface case, though some of the details require extra work.

Note that semi-rational surfaces were also studied in \cite{vS87}.  Our methods are independent of van Straten's.  The end of his thesis contains further information on semi-rational surfaces.  

All varieties are assumed to be quasi-projective over a field of characteristic zero.  In particular, resolution of singularities is true for our varieties.  Macaulayfication is also true for such varieties, but we have not found a way to prove the implication using Kawasaki's construction.  We again note that stating these results in terms of projective closures merely limits our scope to projective varieties; the details will bear out that all results here apply to quasi-projective varieties as well.

\end{section}
\begin{section}{Rational and Semi-Rational Surfaces}
We phrase the main problem for projective varieties.  However, all the techniques carry over to quasi-projective varieties.  In the new language, we are seeking projective closures of quasi-projective varieties with (semi-) rational singularities, the closures having the same types of singularities.

Let $U$ be a demi-normal, quasi-projective surface.  Then, as shown in \cite{Berq14b}, $U$ has a demi-normal projective closure $X$.  Such a surface is Cohen-Macaulay because demi-normal varieties have the $S_2$ property (in addition to being $G_1$ and semi-normal).  Let $p : \overline{X} \rightarrow X$ be the normalization of $X$, and let $\overline{U} = p^{-1}(U)$ be the normalization of $U$, and again note that both $\overline{X}$ and $\overline{U}$ are Cohen-Macaulay because normal varieties are $S_2$ (as well as $R_1$).  

Our first result is that when $U$ is Gorenstein, $U$ has semi-rational singularities if and only if $\overline{U}$ has rational singularities.  We need a preliminary lemma concerning the normalization morphism $p$.

\begin{lemma}  If $p : Z \rightarrow W$ is a finite, surjective morphism of quasi-projective varieties, then the natural morphism $p^*p_*\mathcal{F} \rightarrow \mathcal{F}$ is a surjective morphism of sheaves for any coherent sheaf $\mathcal{F}$.
\proof  This is \cite{Hart87}, III.8.8, together with the fact that for finite, surjective morphisms of quasi-projective varieties, a very ample sheaf of $Z$ over $W$ may be taken as a pullback of a very ample sheaf on $W$ (then apply the projection formula) .  Locally, this is the natural map $M \otimes_A B \rightarrow M$ for a $B$-module $M$.  \qed
\end{lemma}

\begin{proposition} Let $U$ be a demi-normal quasi-projective surface and $p : \overline{U} \rightarrow U$ its normalization.  Suppose $U$ is Gorenstein. Then $U$ has semi-rational singularities if and only if $\overline{U}$ has rational singularities.
\proof  Let $U$ be Gorenstein.  Suppose first that $U$ has semi-rational singularities.  Then $U$ has semi-canonical singularities.  To see this, let $f: V \rightarrow U$ be a resolution.  By definition, we have $f_*\omega_V \cong \omega_U$.  Since $U$ is Gorenstein, $\omega_U$ is invertible.  Thus there is a map $f^*\omega_U \rightarrow \omega_V$, so that there is an effective, $f$-exceptional divisor $E$ such that $\omega_V \cong f^*\omega_U \otimes \mathcal{O}_V(E)$.  This is the definiton of semi-canonical singularities. 

Consider next a commututative diagram $$\begin{CD}
\overline{V}  @>\overline{f}>> \overline{U} \\
@VpVV                          @VqVV \\
V             @>f>>            U \\
\end{CD}$$
Here $f$ is a semi-resolution and $\overline{f}$ the induced resolution of $\overline{U}$.  Note that $V$ is Gorenstein, being semi-smooth.  For the normalization morphism $p$ (and similarly for $q$), we have $\omega_{\overline{V}}(C_0) \cong p^*\omega_V.$  Here $C_0$ is the conductor of $V$ in $\overline{V}$.  Since $U$ has semi-canonical singularities, we may write $$K_V \sim f^*K_U + \sum a_iE_i,$$ where the $E_i$ are $f$-exceptional and the $a_i$ are nonnegative.  

Lifing this relation via $p^*$, and using the above duality statement for $p^*\omega_V$, and transposing, we find that $$K_{\overline{V}} + C_0 \sim \overline{f}^*(K_{\overline{U}} + C) + \sum a_ip^{-1}E_i.$$  In other words, the pair $(\overline{U}, C)$ is log terminal.  This implies that $\overline{U}$ has rational singularities.  See for instance \cite{Kol13}.  Note that the preimages $p^{-1}E_i$ are exactly the $\overline{f}$-exceptional divisors; by definition, a semi-resolution is such that no component of the conductor is exceptional.

For the converse, suppose that $\overline{U}$ has rational singularities.  In terms of the commutative diagram above, this means that $\overline{f}_*\omega_{\overline{V}} \cong \omega_{\overline{U}}.$  We chase $\omega_{\overline{V}}$ in both directions, and we use the duality statement $p_*\omega_{\overline{V}} \cong \mathcal{H}om(p_*\mathcal{O}_{\overline{V}}, \omega_V)$ (and similarly for $q$).  We obtain an isomorphism $$q_*\omega_{\overline{U}} \cong f_*\mathcal{H}om(p_*\mathcal{O}_{\overline{V}}, \omega_V).$$  The right hand side maps to $\mathcal{H}om(q_*\mathcal{O}_{\overline{U}}, f_*\omega_V).$  There is a natural injection $f_*\omega_V \hookrightarrow \omega_U$ obtained by taking the reflexive hull of $f_*\omega_V$.  If we show that this is also surjective, then we will have shown that $U$ has semi-rational singularities.  

What we have shown is that there is a composition, which is the identity (being an identity in codimension 1):  $$\mathcal{H}om(q_*\mathcal{O}_{\overline{U}}, \omega_U) \rightarrow \mathcal{H}om(q_*\mathcal{O}_{\overline{U}}, f_*\omega_V) \rightarrow \mathcal{H}om(q_*\mathcal{O}_{\overline{U}}, \omega_U).$$  In particular, the second map is surjective.  However, each $\mathcal{H}om$ sheaf maps surjectively onto its second factor.  The map is given locally by evaluation at 1.  Thus there is a commutative diagram $$\begin{CD}
\mathcal{H}om(q_*\mathcal{O}_{\overline{U}}, f_*\omega_V) @>>>          \mathcal{H}om(q_*\mathcal{O}_{\overline{U}}, \omega_U) \\
@VVV      @VVV \\                         
f_*\omega_V                                               @>>> \omega_U \\
\end{CD},$$
such that all arrows except possibly the bottom are surjections.  Hence the bottom arrow is a surjection, as required.  We conclude that since $U$ is Cohen-Macaulay (being demi-normal) and $f_*\omega_V \cong \omega_U$, $U$ has semi-rational singularities.\qed
\end{proposition}

We now prove the main result, stated for projective surfaces.  

\begin{theorem}  Suppose that rationalification is true for normal varieties.  Then semi-rationalification is true for Gorenstein demi-normal surfaces.  In other words, if $U$ is a Gorenstein surface with semi-rational sigularities, then $U$ has a projective closure with semi-rational singularities only.
\proof  Let $U$ be a Gorenstein quasi-projective surface with semi-rational singularities, and let $U \hookrightarrow X$ be any projective closure.  Let $p: \overline{X} \rightarrow X$ be the normalization of $X$.  Note that $\overline{U}:= p^{-1}(U) \rightarrow U$ is the normalization of $U$.  By (2.2), $\overline{U}$ has rational singularities.  (As the proof shows, $(\overline{U}, C)$ is log terminal, hence $\overline{U}$ has rational singularities.)  By hypothesis, there exists a proper, birational morphism $f : Y \rightarrow \overline{X}$ such that $Y$ has rational singularities and $f$ is an isomorphism over $\overline{U}$.

Let $\overline{C}$ be the closure in $Y$ of the conductor in $\overline{U}$, and $\overline{D}$ its closed image in $X$.  Note that $\overline{D}$ is the closure in $X$ of the conductor of $U$.  We consider the universal pushout (see \cite{Art70} for its existence) $$\begin{CD}
\overline{C}      @>>> Y \\
@VVV                   @VVV \\
\overline{D}      @>>> X_0 \\
\end{CD}.$$
Here the horizontal arrows are closed immersions and the vertical arrows are finite morphisms.  In fact, $\overline{C} \rightarrow \overline{D}$ is proper and generically finite, and since $\overline{C}$ is a curve, it is proper with finite fibers, hence finite.  In other words, proper with finite fibers implies finite. 

By construction, the universal pushout $g: Y \rightarrow X_0$ is a finite morphism that agrees with $\overline{C} \rightarrow \overline{D}$ and is an isomorphism elsewhere.  Thus in fact $Y$ is the normalization of $X_0$.  Moreover, as $U$ is obtained (by definition) via a pushout $$\begin{CD}
C         @>>> \overline{U} \\
@VVV           @VVV  \\
D         @>>> U \\
\end{CD},$$
and as all these varieties are open subvarieties of the corresponding parts of the pushout $g$, we see that $X_0$ is a projective closure of $U$.  That is, $U$ maps into $X_0$ by the universal property of the pushout, and it does so as an open subvariety.  (It should be checked that $X_0$ is in fact projective.  It is proper, as can be verified using the valuative criterion of properness.  That it has a very ample sheaf follows from the fact that $\overline{C}$, $\overline{D}$, and $Y$ are projective varieties.)

It follows from the main result of \cite{Berq14b} that the projective closure $U \hookrightarrow X_0$ has a demi-normalization that is an isomorphism over $U$.  In other words, there is a finite, birational morphism $X_1 \rightarrow X_0$ that is an isomorphism over $U$ and such that $X_1$ is demi-normal.  In particular, $Y$ is also the normalization of $X_1$.  Now we can conclude that $X_1$ is the desired projective closure of $U$, by (2.2).  This completes the proof.\qed

\end{theorem}

\end{section}
\begin{section}{Criteria for Cohen-Macaulayness}
The proof of Theorem 2.3 carries over directly to varieties of arbitrary dimensions.  In other words, if $X$ and its normalization $\overline{X}$ are both Cohen-Macaulay, then $X$ has semi-rational singularities if and only if $\overline{X}$ has rational singularities.  Thus, a naive approach to proving the main implication would be first to find conditions under which both varieties are Cohen-Macaulay.  We present two such conditions in this section.

Note again that in general there is no connection between $X$ and $\overline{X}$ both being Cohen-Macaulay.  As noted in the introduction, it is possible that $X$ is Cohen-Macaulay even when its normalization is not.  Conversely, taking any surface that is not $S_2$ (such as two planes in 4-space meeting at a point), its normalization is by definition $S_2$ (hence Cohen-Macaulay), so the normalization of $X$ can be Cohen-Macaulay when $X$ is not.  

\begin{proposition}  Consider a demi-normal quasi-projective variety $X$.  Let $p:  \overline{X} \rightarrow X$ be the normalization, and let $p: C \rightarrow D$ be the induced morphism on the conductors.  Suppose that $C$ has rational singularities.  Then if $D$ is normal, $D$ also has rational singularities.  Moreover, under this condition, $X$ is Cohen-Macaulay if and only if so is $\overline{X}$.
\proof  Consider a semi-resolution $f: Y \rightarrow X$ and the induced resolution of singularities $\overline{f}:  \overline{Y} \rightarrow \overline{X}$.  We have a diagram $$\begin{CD}
\overline{Y} @>>> \overline{X} \\
@VqVV             @VpVV \\
Y            @>>> X \\
\end{CD}.$$  Then the morphism of conductors $q: C_0 \rightarrow D_0$ is a finite morphism of smooth varieties; in fact, it is a double cover, ramified along the pinch locus.  By construction of the semi-resolution, $D_0$ maps properly and birationally onto $D$, and hence $C_0$ maps properly and birationally onto $C$.  In other words, when we restrict $f$ and $\overline{f}$ to the conductors, we obtain resolutions of singularities.

Since $q$ is a finite morphism of smooth varieties, $\mathcal{O}_{D_0} \rightarrow q_*\mathcal{O}_{C_0}$ splits locally.  In fact, we can use the trace map to find the splitting (note we are in characteristic zero).  In other words, the trace map is such that the composition $$\mathcal{O}_{D_0} \rightarrow q_*\mathcal{O}_{C_0} \rightarrow \mathcal{O}_{D_0}$$ is the identity.  If we apply $R^if_*$ to this composition, the composition is again the identity:  $$R^if_*\mathcal{O}_{D_0} \rightarrow R^if_*q_*\mathcal{O}_{C_0} \rightarrow R^if_*\mathcal{O}_{D_0}.$$  By the Leray spectral sequence, the middle term is isomorphic to $p_*R^i\overline{f}_*\mathcal{O}_{C_0}$.

By assumption, $C$ has rational singularities.  Thus this last expression, hence the middle term in the above composition, is zero for $i>0$.  It follows that also $R^if_*\mathcal{O}_{D_0}$ is zero for $i> 0$.  In other words, $D$ has rational singularities.  This proves the first assertion.

To complete the proof, we assume that $C$ and $D$ both have rational singularities, and we prove that $X$ has semi-rational singularities if and  only if $\overline{X}$ has rational singularities.  We consider the following commutative diagram of short exact sequences:
$$\begin{CD}
0 @>>> \mathcal{C} @>>>                      \mathcal{O}_Y @>>> \mathcal{O}_{D_0} @>>> 0 \\
@VVV   @VVV                                  @VVV               @VVV                   @VVV \\
0 @>>> q_*\overline{\mathcal{C}} @>>> q_*\mathcal{O}_{\overline{Y}} @>>> q_*\mathcal{O}_{C_0}  @>>> 0  \\
\end{CD}$$
The map on conductor ideal sheaves is an isomorphism.  Moreover, by what we have shown above, the higher direct images of $\mathcal{O}_{D_0}$ and $q_*\mathcal{O}_{C_0}$ under $f$ are all zero.  Therefore, taking the long exact sequences associated to the functor $f_*$ shows that $R^if_*\mathcal{O}_Y = 0$ for $i>1$ if and only if $R^if_*q_*\mathcal{O}_{\overline{Y}} = 0$ for $i>1$.  This last term is isomorphic to $p_*R^i\overline{f}_*\mathcal{O}_{\overline{Y}}$.  Since $p$ is finite, the vanishing of these terms for $i>1$ is equivalent to the vanishing of $R^i\overline{f}_*\mathcal{O}_{\overline{Y}}$ for $i>1$. See (2.1).  Thus our proof is complete, provided we can make the same statement with $R^1f_*$.  However, the splitting of $\mathcal{O}_{D_0} \rightarrow q_*\mathcal{O}_{C_0}$ together with the fact that $\mathcal{C} \rightarrow q_*\overline{\mathcal{C}}$ is an isomorphism implies that the required higher direct images vanish simultaneously.\qed
\end{proposition}

In this case, we were able to show using the definition of (semi-) rational singularities that, provided the conductors have rational singularities, $X$ has semi-rational singularities if and only if $\overline{X}$ has rational singularities. Under weaker conditions on the conductors, we can show that $X$ is Cohen-Macaulay if and only if $\overline{X}$ is Cohen-Macaulay.  This provides another proof that $X$ has semi-rational singularities if and only if $\overline{X}$ has rational singularities.

\begin{proposition}  Suppose, using the definitions in 3.1, that $D$ is normal, and that $C$ and $D$ are both Cohen-Macaulay.  Then $X$ is Cohen-Macaulay if $\overline{X}$ is Cohen-Macaulay.
\proof  We consider the diagram as in 3.1, along with the related diagram of sheaves on $X$:
$$\begin{CD}
0     @>>> \mathcal{C'} @>>>  \mathcal{O}_X  @>>>  \mathcal{O}_D @>>> 0 \\
@VVV       @VVV               @VVV                 @VVV               @VVV \\
0     @>>> p_*\overline{\mathcal{C'}}    @>>>  p_*\mathcal{O}_{\overline{X}} @>>> p_*\mathcal{O}_C @>>> 0 \\
\end{CD}$$
The first nontrivial vertical arrow is again an isomorphism.  We note also that under a finite morphism, the structure sheaf is Cohen-Macaulay if and only if its pushforward is Cohen-Macaulay (as a coherent sheaf of modules).  See \cite{KM98}.  Then taking the long exact sequence in local cohomology, we find that $H^i_x(X, \mathcal{O}_X) = 0$ for $i \leq n-2,$ where $n$ is the dimension of the local ring at a point $x$.  Thus is remains to show that $H^{n-1}_x(X, \mathcal{O}_X) = 0$.  

We have the following commutative diagram of exact sequences:

$$\begin{CD}
0 @>>> H^{n-1}_x(X, \mathcal{C'}) @>>> H^{n-1}_x(X,\mathcal{O}_X)           @>>>            H^{n-1}_x(X, \mathcal{O}_D) \\
@VVV   @VVV                            @VVV          @VVV \\
0 @>>> H^{n-1}_x(X, p_*\overline{\mathcal{C}}) @>>> H^{n-1}_x(X, p_*\mathcal{O}_{\overline{X}})  @>>> H^{n-1}_x(X, p_*\mathcal{O}_C) \\
\end{CD}$$
We note that the middle term in the bottom row is 0, since $\overline{X}$ is Cohen-Macaulay.  Moreover, the final vertical map is an injection.  To see this, it suffices to note that $\mathcal{O}_D \rightarrow p_*\mathcal{O}_C$ splits.  This is true because the corresponding morphism $f_*\mathcal{O}_{D_0} \rightarrow f_*q_*\mathcal{O}_{C_0}$ splits and since $\mathcal{O}_D = f_*\mathcal{O}_{D_0}$ ($D$ is normal and $f: D_0 \rightarrow D$ is a resolution of singularities.)

Now a diagram chase shows that $H^{n-1}_x(X, \mathcal{O}_X) \rightarrow 0 = H^{n-1}_x(X, p_*\mathcal{O}_{\overline{X}})$ is injective.  Thus we have the required vanishing of $H^{n-1}_x(X, \mathcal{O}_X)$.  This finishes the proof.\qed
\end{proposition}

In the final section, we return to the main result.  It will use (3.1).  We note again that the main result is rather incomplete.  Even assuming that $U$ is Gorenstein, what we will show is that, when $U$ has semi-rational singularities, then it has a projective closure $X$ whose normalization has rational singularities.  In the surface case, this was enough to conclude that $X$ itself has semi-rational singularities.  However, as there is no Macaulayfication that preserves the Cohen-Macaulay locus, and as demi-normal varieties of dimension greater than two are not automatically Cohen-Macaulay, it seems that we cannot conclude that $X$ has semi-rational singularities.
\end{section}

\begin{section}{Main Result}  
We restate the version of Koll\'ar's conjecture that we assume.

\noindent
\linebreak
\textbf{Modified Koll\'ar's Conjecture.  }  Suppose that $V$ is quasi-projective, has rational singularities, and that $C$ is a reduced divisor with rational singularities.  Then $V$ has a projective closure $Y$ with rational singularities, and such that the closure of $C$ in $Y$ is Cohen-Macaulay.

\noindent
\linebreak
This seems like a reasonable extension of the original conjecture, which is stated without reference to any divisor $C$.  If fact, what we are stating is similar to the statement that rationalifications of pairs $(X,C)$ exist.  See \cite{Kol13}, 2.82 for the definition of a rational resolution of pairs.  If $X$ has rational singularities in some open set $V$, and $C$ has rational singularities as a subvariety, then we are essentially stating that a rationalification of pairs exists.

\begin{theorem}  Suppose $U$ is a quasi-projective variety with Gorenstein semi-rational singularities.  Then, assuming Koll\'ar's conjecture, $U$ has a projective closure with semi-rational singularities only.
\end{theorem}

Before turning to the proof, we need a preliminary lemma.  When dealing with surfaces, it is imperative that the induced morphism $\overline{C} \rightarrow \overline{D}$ of curves is finite.  All we are given is that it is proper and generically finite.  For curves, this implies that the fibers are finite, hence that $\overline{C} \rightarrow \overline{D}$ is in fact finite.  In general, we will need such a finite morphism on conductors, but finiteness does not come ``for free" the way it does when dealing with curves.  That being said, we have the following:

\begin{lemma}  Suppose that $p: C' \rightarrow D'$ is a proper, generically finite morphism, of generic degree two, between projective varieties. Suppose also that $C'$ is normal.  Let $C \hookrightarrow C'$ and $D \hookrightarrow D'$ be open subvarieties such that $C \rightarrow D$ is finite, and $D$ is normal.  Then there is a finite morphism $C' \rightarrow D''$ to some other normal projective variety $D''$ containing $D$ as an open subvariety, and which agrees with $C \rightarrow D$ when restricted to $C$.
\proof  Since $p$ is a projective morphism, $C'$ is a closed subscheme of some projective space over $D'$:  $$C' \hookrightarrow \textnormal{\textbf{Proj}}_{\mathcal{O}_{D'}}[x_0, \ldots, x_n].$$  We consider the finite morphism $$q: \textnormal{\textbf{Proj}}_{\mathcal{O}_{D'}}[x_0, \ldots, x_n] \rightarrow \textnormal{\textbf{Proj}}_{\mathcal{O}_{D'}}[x_0^2, \ldots, x_n^2].$$   Then we let $D_0$ be the image of $C'$ under $q$.  Locally, $D_0$ is given over $D'=$ Spec $A$ by spectra Spec $A[\{\frac{y_0^2}{y_i^2}\}]$, and by a change of variable (we are in characteristic two), we can assume that the algebra generators are in the fraction field of $A$.  In other words, $D_0$ maps birationally onto $D'$.  We are given a factorization $C' \rightarrow D_0 \rightarrow D'$.  To show that we recover the original morphism $C \rightarrow D$, suppose that $D$ is given locally by Spec $A$ and that $C'$ is given by Spec $B$, where $B$ is a finite $A$-module.  Then $D_0$ is also finite over $D$, as it can be checked that the fibers are finite, and $D_0$ is by construction projective over $D'$.  Being a birational extension, we thus see that $D_0 = D$, since $D$ is assumed to be normal.  Finally, we set $D'' = \overline{D_0}$ to be the normalization of $D_0$.  Then $C'$ still maps onto $D''$, since $C$ is normal, and by the same argument just used, the resulting map agrees with $C \rightarrow D$ when restricted to $C$.\qed
\end{lemma}

With this lemma in hand, we will be able to take a rationalification and glue along its conductor to obtain a variety whose normalization has rational singularities.  The variety will then have semi-rational singularities.  We will need the results of the previous section.  In particular, we will need both conductors (in $X$ and its normalization) to be Cohen-Macaulay, along with the condition that the normalization has rational singularities, to conclude that $X$ is Cohen-Macaulay, and hence has semi-rational singularities.

\noindent
\linebreak
\textbf{Proof of 4.1  }  Suppose that $U$ has Gorenstein semi-rational singularities.  Choose a projective closure $U \hookrightarrow X$, and let $p: \overline{X} \rightarrow X$ be the normalization.  Write $\overline{U} = p^{-1}(U)$ for the normalization of $U$, and let $C \hookrightarrow \overline{U}$ and $D \hookrightarrow U$ be the conductors.  Our first task is to show that $C$ and $D$ both have rational singularities (and in particular, are both Cohen-Macaulay and normal).

As $U$ is Gorenstein, it has an invertible dualizing sheaf.  Moreover, when we pull it back by $p$, we have $p^*\omega_U \cong \omega_{\overline{U}}(C)$.  On the divisor level, this means that $K_{\overline{U}} + C$ is Cartier.  Consider now the diagram $$\begin{CD}
\overline{V} @>\overline{f}>> \overline{U} \\
@VqVV                                        @VpVV \\
V                   @>f>>                 U \\
\end{CD}$$
where $f$ is a semi-resolution of $U$ and $\overline{f}$ the induced resolution of $\overline{U}$.  By construction, the conductor $D_0$ of $V$ is smooth, as is its preimage $C_0$, and $C_0 \rightarrow D_0$ is a double cover, ramified along the pinch locus.   Moreover, $D_0$ maps birationally onto $D$.  It follows that $C_0 \rightarrow C$ is a resolution of singularities.  In addition, $V$ is Gorenstein (having only hypersurface singularities), and when we pull back its dualizing sheaf under $q$, we have $q^*\omega_V \cong \omega_{\overline{V}}(C_0)$.  Again,  on the divisor level, $K_{\overline{V}} + C_0$ is Cartier.

We know that Gorenstein semi-rational singularities are semi-canonical.  In other words, we have a relation $$K_{V} \sim f^*K_U + \sum a_iE_i,$$ where the $E_i$ are $f$-exceptional and the $a_i$ nonnegative.  If we pull back by $q$, we obtain $$K_{\overline{V}} + C_0 \sim \overline{f}^*(K_{\overline{U}} + C) + \sum a_iq^{-1}E_i.$$  In other words, the pair $K_{\overline{U}} + C$ is log terminal.  It follows that $\overline{U}$ has rational singularities.  Moreover, by restricting this last relation to $C_0$ and using the adjunction formula, we find that $C$ itself also has rational singularities.  See \cite{Kol13} for results on adjunction.  Moreover, by \cite{Kol13}, 2.63 and 2.88, we see that $D$ is Cohen-Macaulay because $U$ has semi-canonical singularities.  

It follows that $D$ has rational singularities, since it is Cohen-Macaulay and $C$ has rational singularities.  Too see this, consider the diagram $$\begin{CD}
C_0 @>>> C \\
@VVV          @VVV \\
D_0 @>>> D \\
\end{CD}$$
The horizontal maps are resolutions of singularities, the vertical maps are finite, and $C$ has rational singularities.  As in the proof of the main result for surfaces, this implies that $D$ has rational singularities (we simply chase $\omega_{C_0}$ around the diagram in both directions).  In particular, $D$ is normal.

By hypothesis, there is a rationalification of $\overline{X}$.  In other words, there is a projective morphism $Y \rightarrow \overline{X}$ that is an isomorphism over $\overline{U}$, such that $Y$ has rational singularities, and such that the closure of $C$ in $Y$ is Cohen-Macaulay.  Let this closure be denoted $C'$.  By (4.2), there is a finite morphism $C' \rightarrow D''$ agreeing with $C \rightarrow D$ on $\overline{U}$, and such that $D''$ is normal.  We glue along this morphism:  $$\begin{CD}
C' @>>> Y \\
@VVV      @VgVV \\
D'' @>>> X_0 \\
\end{CD}.$$
The result is a projective variety $X_0$, with normaization $Y$, and which contains $U$ as a dense open subvariety (this last statement is true because $U$ is obtained from $\overline{U}$ by gluing along $C \rightarrow D$).  Since $D''$ is normal, there is a splitting $$\mathcal{O}_{D''} \rightarrow g_*\mathcal{O}_{C'} \rightarrow \mathcal{O}_{D''}.$$  Since $g$ is finite, $g_*\mathcal{O}_{C'}$ is Cohen-Macaulay as an $\mathcal{O}_{D''}$-module.  Then applying local cohomology to this splitting, it is apparent that $D''$ is also Cohen-Macaulay:  $$H^i_x(D'', \mathcal{O}_{D''}) \rightarrow H^i_x(D'', g_*\mathcal{O}_{C'}) = 0 \rightarrow H^i_x(D'', \mathcal{O}_{D''}).$$  For $i>0$, the middle term is zero, and the composition is the identity; therefore the end terms are zero as well.  

It follows from (3.2) that $X_0$ is Cohen-Macaulay.  Finally, since $X_0$ and its normalization are both Cohen-Macaulay, and the normalization has rational singularities, it follows from the usual argument that $X_0$ has semi-rational singularities.  Thus $X_0$ is the desired projective closure of $U$.  This completes the proof.\qed
\end{section}
\nocite{Art70}
\nocite{Kol13}
\nocite{Hart94}
\nocite{Hart87}
\nocite{BH98}
\nocite{GT80}
\nocite{Kov99}
\nocite{LV81}
\nocite{Reid94}
\nocite{dFH09}
\nocite{KSS09}
\nocite{Sch09}
\nocite{vS87}
\nocite{Berq14}
\nocite{GR70}
\nocite{Berq14b}
\nocite{KM98}
\nocite{K00}

\bibliographystyle{plain}
\bibliography{paper}

\end{document}